\documentclass[aop,MSNbibl,number,citesort,dvips]{arximspdf}

\doi{10.1214/11-AOP721} 
\volume{41}
\issue{1}
\pubyear{2013}
\firstpage{109}
\lastpage{133}

\makeatletter
\newcommand{\eqref}[1]{(\ref{#1})}
\newtheorem{theorem}{Theorem}[section]
\newtheorem{corollary}[theorem]{Corollary}
\newtheorem{lemma}[theorem]{Lemma}
\newtheorem{lemmaa}{Lemma}
\newtheorem{proposition}[theorem]{Proposition}
\newproclaim{remark}[theorem]{Remark}
\newproclaim{example}{Example}

\newproclaim{definition}[theorem]{Definition}
\newproclaim{cc}{Continuous case}
\newproclaim{assumption}[theorem]{Assumption}

\makeatother

\begin{document}
\begin{frontmatter}

\title{Functional It{\^o} calculus and stochastic~integral~representation of martingales}
\runtitle{Functional It{\^o} calculus}

\begin{aug}
\author[A]{\fnms{Rama} \snm{Cont}\corref{}\ead[label=e1]{Rama.Cont@upmc.fr}}
\and
\author[B]{\fnms{David-Antoine} \snm{Fourni\'e}\ead[label=e2]{d@vidfournie.com}}
\runauthor{R. Cont and D.-A. Fourni\'e}
\affiliation{CNRS---Universit\'{e} Pierre \& Marie Curie and Columbia University}
\address[A]{Laboratoire de Probabilit\'{e}s et Mod\`{e}les
Ale\'{e}atoires\\
CNRS---Universit\'{e} Pierre \& Marie Curie\\
4 place Jussieu, Case Courrier 188\\
75252 Paris\\
France\\
\printead{e1}}
\address[B]{Department of Mathematics\\
Columbia University\\
New York, New York 10027\\
USA\\
\printead{e2}} 
\end{aug}

\received{\smonth{2} \syear{2010}}
\revised{\smonth{9} \syear{2011}}

%

\begin{abstract}
We develop a nonanticipative calculus for functionals of a
continuous semimartingale, using an extension of the It{\^o} formula to
path-dependent functionals which possess certain directional
derivatives. The construction is based on a pathwise derivative,
introduced by Dupire, for functionals on the space of
right-continuous functions with left limits. We show that this
functional derivative admits a suitable extension to the space of
square-integrable martingales. This extension defines a weak
derivative which is shown to be the inverse of the It{\^o} integral and
which may be viewed as a nonanticipative ``lifting'' of the
Malliavin derivative.

These results lead to a constructive martingale representation
formula for It{\^o} processes. By contrast with the
Clark--Haussmann--Ocone formula, this representation only involves
nonanticipative quantities which may be computed pathwise.
\end{abstract}


\end{frontmatter}

\section{Introduction}\label{sec1}

In the analysis of phenomena with stochastic dynamics, It{\^o}'s
stochastic calculus
\cite{ito44,ito46,dm,meyer76,kw,protter,revuzyor}
has proven to be a powerful and useful tool. A central ingredient of
this calculus is the \textit{It{\^o} formula}~\cite{ito44,ito46,meyer76},
a change of variable formula for functions $f(X_t)$
of a \textit{semimartingale} $X$ which allows one to represent such
quantities in terms of
a stochastic integral.
Given that in many applications such as statistics of processes,
physics or mathematical finance, one is led to consider
path-dependent functionals of a semimartingale
$X$ and its quadratic variation process $[X]$ such as
%
\begin{eqnarray}
\label{examples.eq}
\int_0^t g(t,X_t) \,d[X](t),\qquad
G(t,X_t,[X]_t),\quad \mbox{or}\quad
E[ G(T,X(T),[X](T))|\mathcal{F}_t]\hspace*{-35pt}
\end{eqnarray}
(where $X(t)$ denotes the value at time $t$ and
$X_t=(X(u),u\in[0,t])$ the path up to time $t$), there has been a
sustained interest in extending the framework of stochastic calculus
to such path-dependent functionals.

In this context, the Malliavin calculus
\cite{bismut83,nualart09,malliavin,ocone84,shigekawa80,stroock81,watanabe87}
has proven to be a powerful tool for investigating various
properties of Brownian functionals. Since the construction of
Malliavin derivative does not refer to an underlying {filtration}~$\mathcal{F}_t$,
it naturally leads to representations of functionals
in terms of \textit{anticipative} processes
\cite{clark70,haussmann79,ocone84}. However, in most applications it
is more natural to consider nonanticipative versions of such
representations.

In a recent insightful work, Dupire~\cite{dupire09} has proposed
a method to extend the It{\^o} formula to a functional setting using a pathwise functional derivative
which quantifies the sensitivity of a functional $F_t\dvtx
D([0,t],\mathbb{R})\to\mathbb{R}$ to a variation in the endpoint of
a path $\omega\in D([0,t],\mathbb{R})$.
\[
\nabla_{\omega}F_t(\omega)= \lim_{\varepsilon\to0}\frac{F_t(\omega+
\varepsilon1_{t})- F_t(\omega)}{\varepsilon}.
\]

Building on this insight, we develop hereafter a nonanticipative
calculus~\cite{ContFournie09a} for a class of processes---including
the above examples---which may be represented as
%
\begin{eqnarray}Y(t)= F_t\bigl(
\{X(u),0\leq u\leq t\},\{A(u),0\leq u\leq t\}\bigr) =F_t(X_t,A_t),
\label{Yt.eq}
\end{eqnarray}
where $A$ is the local quadratic variation defined
by $[X](t)=\int_0^t A(u) \,du$, and the functional
\[
F_t\dvtx  D([0,t],\mathbb{R}^d)\times D([0,t],S^+_d)\to
\mathbb{R}
\]
represents the dependence of $Y$ on the path
$X_t=\{X(u),0\leq u\leq t\}$ of $X$ and its quadratic variation.


Our first result (Theorem~\ref{ito.theorem}) is a change of variable
formula for path-dependent functionals of the form \eqref{Yt.eq}.
Introducing $A_t$ as an
additional variable allows us to control the dependence of $Y$ with
respect to the ``quadratic variation'' $[X]$ by requiring smoothness
properties of $F_t$ with respect to the variable $A_t$ in the
supremum norm, without resorting to $p$-variation norms as in
``rough path'' theory~\cite{lyons98}. This allows our result to cover
a wide range of functionals, including the examples in~\eqref{examples.eq}.

We then extend this notion of functional derivative to
\textit{processes}: we show that for $Y$ of the form \eqref{Yt.eq} where
$F$ satisfies some regularity conditions, the process $\nabla_X
Y=\nabla_{\omega}F(X_t,A_t)$ may be defined intrinsically,
independently of the choice of~$F$ in \eqref{Yt.eq}. The operator
$\nabla_X$ is shown to admit an extension to the space of
square-integrable martingales, which is the inverse of the It{\^o}
integral with respect to $X$: for $\phi\in\mathcal{L}^2(X),
\nabla_X( \int\phi\cdot dX ) = \phi$ (Theorem
\ref{weakderivative.theorem}). In particular, we obtain a
constructive version of the martingale representation theorem
(Theorem~\ref{L2martingalerepresentation.theorem}), which states
that for any square-integrable $\mathcal{F}^X_t$-martingale $Y$,
\[
Y(T) = Y(0) + \int_0^T \nabla_XY \cdot dX,\qquad \mathbb{P}\mbox{-a.s.}
\]
This formula can be seen as a nonanticipative counterpart of the
Clark--Hauss\-mann--Ocone formula
\cite{clark70,haussmann78,haussmann79,karatzasocone91,ocone84}. The
integrand $\nabla_XY$ is an adapted process which may be computed
pathwise, so this formula is more amenable to numerical
computations than those based on Malliavin calculus.

Finally, we show that this functional derivative $\nabla_X$ may
be viewed as a nonanticipative ``lifting'' of the Malliavin
derivative (Theorem~\ref{relevement.theorem}): for
square-integrable martingales $Y$ whose terminal values is
differentiable in the sense of Malliavin $Y(T)\in\mathbf{D}^{1,2}$, we
show that $\nabla_XY(t)=E[\mathbb{D}_tH|\mathcal{F}_t]$.

These results provide a rigorous mathematical framework for
developing and extending the ideas proposed by Dupire
\cite{dupire09} for a large class of functionals.
In particular, unlike the results derived from the pathwise approach
presented in~\cite{ContFournie09c,dupire09},
Theorems
\ref{weakderivative.theorem} and \ref
{L2martingalerepresentation.theorem} do not require any pathwise
regularity of the functionals and hold for nonanticipative
square-integrable processes, including stochastic integrals and
functionals which may depend on the
quadratic variation of the process.

\section{Functional representation of nonanticipative processes}
\label{functionalrepresentation.sec}

Let $X\dvtx [0,T]\times\Omega\mapsto\mathbb{R}^d$ be a continuous,
$\mathbb{R}^d$-valued semimartingale defined on a filtered
probability space $(\Omega,\mathcal{F},\mathcal{F}_t,\mathbb{P})$ assumed
to satisfy the usual hypotheses~\cite{dm}. Denote by~$\mathcal{P}$
(resp., $\mathcal{O}$) the associated \textit{predictable} (resp., optional)
sigma-algebra on $[0,T]$. $\mathcal{ F}^X_{t}$ denotes the
($\mathbb{P}$-completed) natural filtration of $X$.
The paths of
$X$ then lie in $C_0([0,T],\mathbb{R}^d)$, which we will view as a
subspace of $D([0,T],\mathbb{R}^d)$ the space of cadlag functions
with values in $\mathbb{R}^d$. We denote by
$[X]=([X^i,X^j],i,j=1,\ldots,d)$ the quadratic (co-)variation process
associated to $X$,
taking values in the set $S^+_d$ of
positive $d \times d$ matrices. We assume that
%
\begin{equation}[X](t)=\int_0^t
A(s)\,ds \label{quadraticvariationrepresentation.eq}
\end{equation}
for
some cadlag process $A$ with values in $S^+_d$. Note that $A$ need
not be a semimartingale. The paths of $A$ lie in
$\mathcal{S}_t=D([0,t], S^+_d)$, the space of cadlag functions with
values $S^+_d$.
\subsection{Horizontal extension and vertical perturbation of a
path}\label{extensions.sec}

Consider a path $x\in D([0,T]),\mathbb{R}^d)$ and denote by
$x_t=(x(u), 0\leq u\leq t)\in D([0,t],\mathbb{R}^d)$ its restriction
to $[0,t]$ for $t<T$. For a process $X$ we shall similarly denote
$X(t)$ its value at $t$ and
$X_t=(X(u), 0\leq u\leq t)$ its path on $[0,t]$.

For $h\geq0$, we define the \textit{horizontal} extension
$x_{t,h}\in D([0,t+h],\mathbb{R}^d)$ of $x_t$ to $[0,t+h]$ as
%
\begin{equation}
x_{t,h}(u)=x(u),\qquad u\in[0,t] ;\qquad x_{t,h}(u)=x(t),\qquad
u\in\,]t,t+h].\vadjust{\goodbreak}
\end{equation}
For $h\in\mathbb{R}^d$, we define the \textit{vertical}
perturbation $x^h_t$ of $x_t$ as the cadlag path obtained by
shifting the endpoint by $h$.
%
\begin{eqnarray}x^h_t(u)=x_t(u),\qquad u\in[0,t[,\qquad x^h_t(t)=x(t)+h,
\end{eqnarray}
or, in other words, $x^h_t(u)=x_t(u)+h 1_{t=u}$.

\subsection{Adapted processes as nonanticipative functionals}
A process $Y\dvtx\break [0,T]\times\Omega\mapsto\mathbb{R}^d$ adapted to
$\mathcal{ F}^X_t$ may be represented as
%
\begin{equation}Y(t)= F_t\bigl(
\{X(u),0\leq u\leq t\},\{A(u),0\leq u\leq t\}\bigr)
=F_t(X_t,A_t),\label{representation.eq}
\end{equation}
where
$F=(F_t)_{t\in[0,T]}$ is a family of functionals
\[F_t\dvtx
D([0,t],\mathbb{R}^d)\times\mathcal{S}_t\to\mathbb{R}
\]
representing the dependence of $Y(t)$ on the underlying path of $X$
and its quadratic variation.

Since $Y$ is nonanticipative, $Y(t,\omega)$ only depends on the
restriction $\omega_t$ of $\omega$ on $[0,t]$. This motivates the
following definition:
\begin{definition}[(Nonanticipative functional)]\label{functional.def} A
nonanticipative functional is a
family of functionals $F=(F_t)_{t\in[0,T]}$ where
\begin{eqnarray*}F_t \dvtx
D([0,t],\mathbb{R}^d)\times D([0,t],{S}^+_d)&\mapsto&
\mathbb{R},\\
(x,v)& \to& F_t(x,v)
\end{eqnarray*}
is measurable with respect to
$\mathcal{B}_t,$ the canonical filtration on
$D([0,t],\mathbb{R}^d)\times D([0,t],{S}^+_d)$.
\end{definition}

We can also view $F=(F_t)_{t\in[0,T]}$ as a map defined on the
space $\Upsilon$ of \textit{stopped paths}
%
\begin{equation}\Upsilon= \{
(t,\omega_{t,T-t}), (t,\omega)\in[0,T]\times
D([0,T],\mathbb{R}^d\times{S}^+_d ) \}.
\end{equation}
Whenever the context is
clear, we will denote a generic element $(t,\omega)\in\Upsilon$
simply by its second component, the path $\omega$ stopped at $t$.
$\Upsilon$ can also be identified with the ``vector bundle''
%
\begin{equation}
\Lambda= \bigcup_{t\in[0,T]}D([0,t],\mathbb{R}^d) \times
D([0,t],{S}^+_d).
\end{equation}

A natural distance on the space $\Upsilon$ of stopped paths is given
by
%
\begin{equation}d_\infty( (t,\omega), (t',\omega'))= |t-t'|+ \sup
_{u\in
[0,T]} |\omega_{t,T-t}(u)-\omega'_{t',T-t'}(u)|.
\end{equation}
$(\Upsilon,d_\infty)$ is then a metric space, a closed subspace of
$([0,T]\times D([0,T],\mathbb{R}^d\times{S}^+_d ),\|\cdot\|_\infty)$ for
the product topology.

Introducing the process $A$ as an additional variable
may seem redundant at this stage: indeed $A(t)$ is itself $\mathcal{F}_t$-measurable,
that is, a functional of $X_t$. However, it is not a
\textit{continuous} functional on $(\Upsilon,d_\infty)$. Introducing
$A_t$ as a second argument in the functional will allow us to
control the regularity of $Y$ with respect to $[X]_t=\int_0^t A(u)
\,du$
simply by requiring
continuity of $F_t$ in supremum or $L^p$ norms with respect to the
``lifted process'' $(X,A)$; see Section~\ref{regularity.sec}. This
idea is analogous in some ways to the approach of rough path theory
\cite{lyons98}, although here we do not resort to $p$-variation norms.

If $Y$ is a $\mathcal{B}_t$-{predictable} process, then \cite[Volume
I, paragraph 97]{dm}
\[
\forall t\in[0,T],\qquad  Y(t,\omega)=Y(t,\omega_{t-}),
\]
where
$\omega_{t-}$ denotes the path defined on $[0,t]$ by
\[
\omega_{t-}(u)=\omega(u),\qquad u\in[0,t[,\qquad \omega_{t-}(t)=\omega(t-).
\]
Note that $\omega_{t-}$ is cadlag and should \textit{not} be
confused with the caglad path $u\mapsto\omega(u-)$.

The functionals discussed in the introduction depend on the process $A$
via $[X]=\int_0^{\cdot} A(t) \,dt$. In particular, they satisfy the condition
$F_t(X_t,A_t)=F_t(X_t,A_{t-})$.
Accordingly, we will assume throughout the paper that all functionals
$ F_t\dvtx  D([0,t], \mathbb{R}^d) \times
\mathcal{S}_t\to\mathbb{R}\nonumber$ considered have
``predictable'' dependence with respect to the second argument,
%
\begin{equation}
\quad\forall t \in[0,T], \forall(x,v) \in D([0,t],\mathbb{R}^d)
\times\mathcal{S}_t,\qquad
F_t(x_t,v_t)=F_t(x_t,v_{t-}).\label{predictable.eq}
\end{equation}

\subsection{Continuity for nonanticipative functionals}\label{regularity.sec}
We now define a notion of (left) continuity for nonanticipative
functionals.
\begin{definition}[(Continuity at fixed times)]\label{fixedtimecontinuous.def}
A functional $F$ defined on $\Upsilon$ is said to be continuous at
fixed times for
the $d_\infty$ metric if and only if
%
\begin{eqnarray}
&&\forall t \in[0,T), \forall\varepsilon> 0,
\forall(x,v)\in
D([0,t],\mathbb{R}^d) \times\mathcal{S}_t,\nonumber\\
&&\quad\exists\eta> 0,
(x',v') \in D([0,t],\mathbb{R}^d) \times
\mathcal{S}_{t},
\\
 &&\qquad d_\infty((x,v),(x',v')) < \eta
\quad\Rightarrow\quad
|F_t(x,v)-F_{t}(x',v')| < \varepsilon.\nonumber
\end{eqnarray}
\end{definition}

We now define a notion of joint continuity with respect to time and
the underlying path:
\begin{definition}[(Continuous functionals)]
A nonanticipative functional $F=(F_t)_{t\in[0,T)}$ is said to be
continuous at $(x,v)\in D([0,t],\mathbb{R}^d)\times\mathcal{S}_t$
if
%
\begin{eqnarray}&&\forall\varepsilon> 0, \exists\eta> 0, \forall(x',v')
\in
\Upsilon,
\nonumber
\\[-8pt]
\\[-8pt]
\nonumber
&& \qquad d_\infty((x,v),(x',v')) < \eta\quad\Rightarrow\quad
|F_t(x,v)-F_{t'}(x',v')| < \varepsilon.
\end{eqnarray}
We denote by
$\mathbb{C}^{0,0}([0,T) )$ the set of continuous nonanticipative functionals
on $\Upsilon$.
\end{definition}
\begin{definition}[(Left-continuous functionals)]\label{lcontinuous.def}
A nonanticipative functional $F=(F_t,t\in[0,T) )$ is said to be
\mbox{left-continuous} if for each $t\in[0,T),$ $F_t\dvtx
D([0,t], \mathbb{R}^d) \times\mathcal{S}_t\to\mathbb{R}$ in the sup
norm and
%
\begin{eqnarray}
&&\forall\varepsilon> 0, \forall(x,v)\in
D([0,t],\mathbb{R}^d) \times\mathcal{S}_t,\nonumber \\
&&\quad\exists\eta> 0,
\forall h \in[0,t], \forall(x',v') \in
D([0,t-h],\mathbb{R}^d) \times\mathcal{S}_{t-h},
\\
&&\qquad d_\infty((x,v),(x',v')) < \eta\quad\Rightarrow\quad|F_t(x,v)-F_{t-h}(x',v')|
< \varepsilon.\nonumber
\end{eqnarray}
We denote by
$\mathbb{C}^{0,0}_l([0,T) )$ the set of left-continuous
functionals.
\end{definition}

We define analogously the class of right-continuous functionals
$\mathbb{C}^{0,0}_r([0,T) )$.

We call a functional ``boundedness preserving'' if it is bounded
on each bounded set of paths:
\begin{definition}[(Boundedness-preserving functionals)]\label{boundedness.def}
Define $\mathbb{B}([0,T))$ as the set of nonanticipative
functionals $F$ such that for every compact subset $K$ of~$\mathbb{R}^d$, every $R >0$ and $t_0 < T$,
%
\begin{eqnarray}\label{boundedpreserving.eq}
&&\exists
C_{K,R,t_0}>0, \forall t \leq t_0, \forall(x,v) \in
D([0,t],K) \times\mathcal{S}_t,
\nonumber
\\[-8pt]
\\[-8pt]
\nonumber
&&\qquad\sup_{s \in[0,t]} |v(s)| < R
\quad\Rightarrow\quad|F_t(x,v)| < C_{K,R,t_0}.
\end{eqnarray}
\end{definition}
%
\subsection{Measurability properties}
Composing a nonanticipative functional $F$ with the process $(X,A)$
yields an $\mathcal{F}_t$-adapted process $Y(t)=F_t(X_t,A_t)$. The
results below link the measurability and pathwise regularity of $Y$
to the regularity of the functional $F$.
\begin{lemma}[(Pathwise regularity)] If $F\in\mathbb{C}^{0,0}_l$, then
for any $(x,v)\in
D([0,T], \mathbb{R}^d) \times\mathcal{S}_T$, the path $t\mapsto
F_t(x_{t-},v_{t-})$ is left-continuous.
\label{cadlag.prop}
\end{lemma}
\begin{pf}
Let $F\in\mathbb{C}^{0,0}_l$ and $t \in[0,T)$. For $h> 0$
sufficiently small,
%
%
\begin{eqnarray}
d_{\infty}((x_{t-h},v_{t-h}),(x_{t-},v_{t-}))&=&\sup_{u \in(t-h,t)}
|x(u)-x(t-h)|
\nonumber
\\[-8pt]
\\[-8pt]
\nonumber
&&{} + \sup_{u \in(t-h,t)} |v(u)-v(t-h)| + h.
\end{eqnarray}
Since $x$ and $v$ are cadlag, this quantity converges to 0 as $h
\rightarrow0+$, so
\[
F_{t-h}(x_{t-h},v_{t-h})-F_t(x_{t-},v_{t-}) \mathop{\to}^{h\to0^+} 0,
\]
so $t\mapsto
F_t(x_{t-},v_{t-})$ is left-continuous.
\end{pf}

\begin{theorem} \label{Measurability}
\textup{(i)} If $F$ is continuous at fixed times, then the process $Y$
defined by $Y((x,v),t)=F_t(x_{t},v_{t})$ is adapted.\vspace*{-6pt}
\begin{longlist}[(iii)]
\item[(ii)] If $F\in\mathbb{C}^{0,0}_l([0,T) )$, then the process
$Z(t)=F_t(X_t,A_t)$ is optional.
\item[(iii)] If $F\in\mathbb{C}^{0,0}_l([0,T) )$, and if either $A$ is
continuous or $F$ verifies \eqref{predictable.eq}, then~$Z$ is a
predictable process.
\end{longlist}
\end{theorem}

In particular, any $F\in\mathbb{C}^{0,0}_l$ is a nonanticipative
functional in the sense of Definition~\ref{functional.def}. We
propose an easy-to-read proof of points (i) and (iii) in the case
where~$A$ is continuous. The (more technical) proof for the cadlag
case is given in the \hyperref[measurabilityproof.sec]{Appendix}.
\begin{cc*}
Assume that $F$ is continuous at fixed times and that the paths of
$(X,A)$ are almost-surely continuous. Let us prove that $Y$ is $\mathcal{F}_t$-adapted:
$X({t})$ is $\mathcal{F}_t$-measurable. Introduce the partition
$t^{i}_{n}=\frac{iT}{2^{n}},i=0,\ldots,2^{n}$ of $[0,T]$, as well as the
following piecewise-constant approximations of $X$ and~$A$:
%
%
\begin{eqnarray}
X^{n}(t)&=&\sum_{k=0}^{2^{n}} X(t^{n}_{k})1_{[t^{n}_{k},t^{n}_{k+1})}(t)+
X_{T}1_{\{T\}}(t),
\nonumber
\\[-8pt]
\\[-8pt]
\nonumber
A^{n}(t)&=&\sum_{k=0}^{2^{n}} A(t^{n}_{k})1_{[t^{n}_{k},t^{n}_{k+1})}(t)+
A_{T}1_{\{T\}}(t).
\end{eqnarray}
The random variable $Y^{n}(t)=F_{t}(X^{n}_{t},A^{n}_{t})$ is a
continuous function of the random variables
$\{X(t^{n}_{k}),A(t^{n}_{k}),t^{n}_{k}\leq t\}$ and hence is $\mathcal
{F}_t$-measurable. The representation above shows in fact
that $Y^{n}(t)$ is $\mathcal{F}_t$-measurable.
$X^{n}_{t}$ and $A^{n}_{t}$ converge respectively to $X_{t}$ and
$A_{t}$ almost-surely
so $ Y^{n}(t) \mathop{\to}^{n\to\infty} Y(t)$ a.s., and hence
$Y(t)$ is
$\mathcal{F}_t$-measurable.

(i) implies point (iii) since the path of $Z$ are left-continuous by
Lemma~\ref{cadlag.prop}.
\end{cc*}

\section{Pathwise derivatives of nonanticipative functionals}
\label{derivative.sec}
\subsection{Horizontal and vertical derivatives}
We now define pathwise derivatives for a nonanticipative functional, following
Dupire~\cite{dupire09}.
\begin{definition}[(Horizontal derivative)]\label{horizontalderivative.def}
The \textit{horizontal derivative} at $(x,v)\in
D([0,t],\mathbb{R}^d)\times\mathcal{S}_t$ of nonanticipative
functional $F=(F_t)_{t\in[0,T) }$ is defined as
%
\begin{equation}\mathcal{D}_tF(x,v) = \lim_{h\to0^+} \frac{F_{t+h}(
x_{t,h},v_{t,h})-F_t(x_t,v_t)}{h}\label{horizontalderivative.eq}
\end{equation}
if the corresponding limit exists. If
\eqref{horizontalderivative.eq} is defined for all
$(x,v)\in\Upsilon$, the map
%
\begin{eqnarray}\mathcal{D}_tF \dvtx
D([0,t],\mathbb{R}^d)\times
\mathcal{S}_t &\mapsto& \mathbb{R}^d,
\nonumber
\\[-8pt]
\\[-8pt]
\nonumber
(x,v)&\to& \mathcal{D}_tF(x,v)
\end{eqnarray}
defines a nonanticipative
functional $\mathcal{D}F=(\mathcal{D}_tF)_{t\in[0,T]}$, the \textit{
horizontal derivative} of $F$.\vadjust{\goodbreak}
\end{definition}

Note that our definition \eqref{horizontalderivative.eq} is different
from the one
in~\cite{dupire09} where the case $F(x,v)=G(x)$ is considered.


Dupire~\cite{dupire09} also introduced a pathwise spatial derivative
for such functionals, which we now introduce. Denote $(e_i,i=1,\ldots,d)$
the canonical basis in~$\mathbb{R}^d$.
\begin{definition} A nonanticipative
functional $F=(F_t)_{t\in[0,T) }$ is said to be
\textit{vertically differentiable} at $(x,v)\in
D([0,t]),\mathbb{R}^d)\times D([0,t],S^+_d)$ if
%
\begin{eqnarray*}\mathbb{R}^d&\mapsto&\mathbb{R},
\\
e&\to& F_t(x^{e}_t,v_t)
\end{eqnarray*}
is differentiable at $0$.
Its gradient at $0$,
%
\begin{eqnarray} \label{verticalderivative.eq}
{\nabla}_xF_t (x, v) = \bigl(\partial_iF_t(x,v), i=1,\ldots,d\bigr)
\nonumber
\\[-8pt]
\\[-8pt]
\eqntext{\mbox{where }\partial_iF_t(x,v)=\displaystyle\lim_{h\to0}\frac
{F_t(x^{he_i}_t,v)-F_t(x,v)}{h}}
\end{eqnarray}
is called the \textit{vertical derivative} of $F_t$ at $(x,v)$. If
\eqref{verticalderivative.eq} is defined for all $(x,v)\in\Upsilon$,
the maps
%
\begin{eqnarray}
&\displaystyle{\nabla}_xF \dvtx  D([0,t],\mathbb{R}^d)\times
\mathcal{S}_t \mapsto \mathbb{R}^d,&
\nonumber
\\[-8pt]
\\[-8pt]
\nonumber
&\displaystyle(x,v)\to {\nabla}_xF_t(x,v)&
\end{eqnarray}
define a nonanticipative
functional ${\nabla}_xF=({\nabla}_xF_t)_{t\in[0,T]}$, the \textit{
vertical derivative} of $F$. $F$ is then said to be \textit{vertically
differentiable} on $\Upsilon$.
\label{verticalderivative.def}
\end{definition}
\begin{remark}$\partial_iF_t(x,v)$ is simply the directional derivative
of $F_t$
in the direction $(1_{\{t\}} e_i,0)$. Note that this involves evaluating $F$ at
cadlag perturbations of the path $x$, even if $x$ is continuous.
\end{remark}
\begin{remark}
If $F_t(x,v)=f(t,x(t))$ with $f\in C^{1,1}([0,T) \times
\mathbb{R}^d)$, then we retrieve the usual partial derivatives
\[
{\mathcal D}_tF(x,v)= \partial_tf(t,X(t)),\qquad \nabla_xF_t(X_t,A_t)=
\nabla_xf(t,X(t)).
\]
\end{remark}
\begin{remark}
Bismut~\cite{bismut83} considered directional derivatives of
functionals on $D([0,T],\mathbb{R}^d)$ in the direction of purely
discontinuous (e.g., piecewise constant) functions with finite
variation, which is similar to Definition~\ref{verticalderivative.def}.
This notion, used in~\cite{bismut83} to derive an integration by
parts formula for pure-jump processes, is natural in the context of
discontinuous semimartingales.
We will show that the directional derivative \eqref
{verticalderivative.eq} also
intervenes naturally when the underlying process $X$ is
\textit{continuous}, which is less obvious.
\end{remark}
\begin{definition}[(Regular functionals)] \label{Cab.def}
Define $\mathbb{C}^{1,k}([0,T))$ as the set of functionals $F \in
\mathbb{C}^{0,0}_l$ which are:\vadjust{\goodbreak}
\begin{itemize}
\item horizontally differentiable with $\mathcal{D}_tF$ continuous at fixed
times;
\item$k$ times vertically differentiable with $\nabla^j_xF \in
\mathbb{C}^{0,0}_l([0,T))$ for $j=1,\ldots,k$.
\end{itemize}
Define $\mathbb{C}^{1,k}_b([0,T) )$ as the set of functionals $F \in
\mathbb{C}^{1,2}$ such that $\mathcal{D}F,\nabla_xF,\ldots, \break\nabla^k_xF
\in\mathbb{B}([0,T))$.
\end{definition}

We denote $\mathbb{C}^{1,\infty}([0,T) )=\bigcap_{k\geq1}
\mathbb{C}^{1,k}([0,T) $.

Note that this notion of regularity only involves directional
derivatives with
respect to \textit{local} perturbations of paths, so $\nabla_x F$ and
$\mathcal{D}_tF$ seems to contain \textit{less} information on the
behavior of $F$ than, say, the Fr\'echet derivative which considers
perturbations in all directions in $C_0([0,T],\mathbb{R}^d)$ or the
Malliavin derivative~\cite{malliavin78,malliavin} which examines
perturbations in the direction of all absolutely continuous
functions. Nevertheless we will show in Section
\ref{functionalito.sec} that knowledge of
$\mathcal{D}F,\nabla_xF,\nabla^2_xF$ along the paths of $X$
is sufficient to reconstitute the path of
$Y(t)=F_t(X_t,A_t)$.
\begin{example}[(Smooth functions)]
In the case where $F$ reduces to a
smooth \textit{function} of $X(t)$,
%
\begin{equation}F_t(x_t,v_t)=f(t,x(t)),
\end{equation}
where $f \in C^{1,k}([0,T]\times\mathbb{R}^d)$, the pathwise
derivatives reduce to the usual ones
%
\begin{eqnarray}\mathcal{D}_tF(x_t,v_t)=\partial_tf(t,x(t)),\qquad
\nabla^j_xF_t(x_t,v_t)=\partial^j_xf(t,x(t)).
\end{eqnarray}
In fact to have $F \in\mathbb{C}^{1,k}$ we just need $f$ to be
right-differentiable in the time variable, with right-derivative
$\partial_tf(t,\cdot)$ which is continuous in the space variable and~$f$,
$\nabla f$ and $\nabla^2f$ to be jointly left-continuous in $t$
and continuous in the space variable.
\end{example}
\begin{example}[(Cylindrical functionals)]\label{example.cylindrical} Let
$g \in
C^0(\mathbb{R}^d,\mathbb{R}), h\in C^k(\mathbb{R}^d,\mathbb{R} )$
with $h(0)=0$. Then
\[
F_t(\omega)=
h\bigl(\omega(t)-\omega(t_n-)\bigr) 1_{t\geq t_n}
g(\omega(t_1-),\omega(t_2-),\ldots,\omega(t_n-))
\]
is in $
\mathbb{C}^{1,k}_b$ with $\mathcal{D}_tF(\omega)=0$ and
\begin{eqnarray*}
&&\forall j=1,\ldots,k,\\
&&\qquad \nabla^j_{\omega}F_t(\omega)=
h^{(j)}\bigl(\omega(t)-\omega(t_n-)\bigr) 1_{t\geq t_n}
g(\omega(t_1-),\omega(t_2-)\ldots,\omega(t_n-)).
\end{eqnarray*}
\end{example}
\begin{example}[(Integrals with respect to quadratic variation)]\label{integral.ex}
A process $Y(t)=\int_0^t g(X(u)) \,d[X](u)$ where $g \in
C^{0}(\mathbb{R}^d)$ may be represented by the functional
%
\begin{equation}F_t(x_t,v_t)=\int_0^t g(x(u)) v(u)\,du.
\end{equation}
It is readily observed that $F \in\mathbb{C}^{1,\infty}_b$, with
%
\begin{equation}\mathcal{D}_tF(x_t,v_t)=g(x(t))v(t),\qquad \nabla
^j_xF_t(x_t,v_t)=0.
\end{equation}
\end{example}
\begin{example}\label{quadratic.ex} The martingale $Y(t)=X(t)^2-[X](t)$ is represented by
the functional
%
\begin{equation}
F_t(x_t,v_t)=x(t)^2-\int_0^t v(u)\,du.
\end{equation}
Then $F \in
\mathbb{C}^{1,\infty}_b$ with
%
\begin{eqnarray}
\mathcal{D}_tF(x,v)&=&-v(t),\qquad
\nabla_xF_t(x_t,v_t)=2x(t),
\nonumber
\\[-8pt]
\\[-8pt]
\nonumber
\nabla^2_xF_t(x_t,v_t)&=&2,\qquad \nabla^j_xF_t(x_t,v_t)=0,\qquad j\geq3.
\end{eqnarray}
\end{example}
\begin{example}\label{doleans.ex} $Y=\exp(X-[X]/2)$ may be represented as $Y(t)=F(X_t)$
%
\begin{equation}
F_t(x_t,v_t)=e^{x(t)-{1}/{2}\int_0^tv(u)\,du}.
\end{equation}
Elementary
computations show that $F \in\mathbb{C}^{1,\infty}_b$ with
%
\begin{equation}
\mathcal{D}_tF(x,v)=-\frac{1}{2}v(t)F_t(x,v),\qquad
\nabla^j_xF_t(x_t,v_t)=F_t(x_t,v_t).
\end{equation}
\end{example}

Note that, although $A_t$ may be expressed as a functional of $X_t$,
this functional is not continuous
and without introducing the second variable $v\in\mathcal{S}_t$, it is
not possible to represent Examples~\ref{integral.ex},
\ref{quadratic.ex} and~\ref{doleans.ex} as a left-continuous functional
of $x$ alone.

\subsection{Obstructions to regularity}
It is instructive to observe what prevents a functional from being
regular in the sense of Definition~\ref{Cab.def}. The examples below
illustrate the fundamental obstructions to regularity:
\begin{example}[(Delayed functionals)] Let $\varepsilon>0$.
$F_t(x_t,v_t)=x(t-\varepsilon)$ defines a $\mathbb{C}^{0,\infty}_b$
functional. All vertical derivatives are 0. However, $F$ fails to be
horizontally differentiable.
\end{example}

\begin{example}[(Jump of $x$ at the current time)]
$F_t(x_t,v_t)=x(t)-x(t-)$ defines a functional which is infinitely
differentiable and has regular pathwise derivatives
%
\begin{eqnarray}
\mathcal{D}_tF(x_t,v_t)=0,\qquad \nabla_xF_t(x_t,v_t)=1.
\end{eqnarray}
However,
the functional itself fails to be $\mathbb{C}^{0,0}_l$.
\end{example}
\begin{example}[(Jump of $x$ at a fixed time)]
$F_t(x_t,v_t)=1_{t \geq t_0}(x(t_0)-x(t_0-))$ defines a
functional in $\mathbb{C}^{0,0}_l$ which admits horizontal
and vertical derivatives at any order at each point $(x,v)$.
However, $\nabla_xF_t(x_t,v_t)=1_{t=t_0}$ fails to be either right- or
left-continuous, so $F$ is not $\mathbb{C}^{0,1}$ in the sense of
Definition~\ref{verticalderivative.def}.\vadjust{\goodbreak} 
\end{example}
\begin{example}[(Maximum)]
$F_t(x_t,v_t)=\sup_{s \leq t} x(s)$ is $\mathbb{C}^{0,0}_l$ but
fails to be vertically differentiable on the set
\[
\Bigl\{ (x_t,v_t)\in D([0,t],\mathbb{R}^d)\times\mathcal{S}_t,
x(t)=\sup_{s \leq t} x(s)\Bigr\}.
\]
\end{example}

\section{Functional It{\^o} calculus}\label{functionalito.sec}

\subsection{Functional It{\^o} formula}
We are now ready to prove our first main result, which is a change
of variable formula for nonanticipative functionals of a
semimartingale~\cite{ContFournie09a,dupire09}:
\begin{theorem}\label{ito.theorem}
For any nonanticipative functional $F \in\mathbb{C}^{1,2}_b$
verifying \eqref{predictable.eq} and any $t \in[0, T) $,
%
%
\begin{eqnarray}\label{functional.ito.eq}
&&F_t(X_t,A_t)-F_0(X_0,A_0) \nonumber\\
&&\qquad=\int_{0}^{t} \mathcal{D}_uF(X_u,A_u) \,du+
\int_{0}^{t}
{\nabla}_xF_u(X_u,A_u)\cdot dX(u)
\\
&&\qquad\quad{}+ \int_0^t \frac{1}{2} \operatorname{tr}(\nabla^2_xF_u(X_u,A_u) \,d[X](u))\qquad
a.s.\nonumber
\end{eqnarray}
In particular, for any $F\in\mathbb{C}^{1,2}_b$, $Y(t)=F_t(X_t,A_t)$
is a
semimartingale.
\end{theorem}

Theorem~\ref{ito.theorem} shows that, for a regular
functional $F\in
\mathbb{C}^{1,2}([0,T))$, the process $Y=F(X,A)$ may be
reconstructed from the second-order jet $(\mathcal{D}F,\nabla_xF,\break\nabla_x^2F)$ of $F$ along the paths of $X$.

\begin{pf}
Let us first assume that $X$ does not exit a compact set $K$ and that
$\|A\|_\infty\leq R$ for some $R> 0$. Let us introduce a sequence of
random partitions $(\tau^{n}_{k},k=0,\ldots,k(n))$ of $[0,t]$, by adding the
jump times of $A$ to the dyadic partition $(t^{n}_{i}=\frac
{it}{2^{n}},i=0,\ldots,2^{n})$,
%
%
\begin{equation}
\qquad\tau^{n}_{0}=0,\qquad
\tau^{n}_{k}=\inf\biggl\{s > \tau^{n}_{k-1} | 2^{n}s \in\mathbb{N} \mbox{ or } |A(s)-A(s-)| > \frac{1}{n} \biggr\} \wedge t.
\end{equation}
The following arguments apply pathwise. Lemma~\ref{Approximation}
ensures that
\[
\eta_n=\sup\biggl\{|A(u)-A(\tau^n_i)|+|X(u)-X(\tau^n_i)|+\frac
{t}{2^n},i\leq
2^n,u\in[\tau^n_{i},\tau^n_{i+1})\biggr\}\mathop{\rightarrow
}_{n\rightarrow
\infty} 0.
\]
Denote
$_{n}X =\sum_{i=0}^{\infty}X(\tau^n_{i+1})1_{[\tau^n_{i},\tau
^n_{i+1})} +
X(t)1_{\{t\}}$
which is a cadlag piecewise constant approximation of $X_t$, and
$_nA=\sum_{i=0}^{\infty}A(\tau^n_i)1_{[\tau^n_{i},\tau^n_{i+1})} +
A(t)1_{\{t\}}$ which is an adapted cadlag piecewise constant\vadjust{\goodbreak}
approximation of $A_t$.
Denote $h^n_i=\tau^{n}_{i+1}-\tau_i^n$. Start with the decomposition
%
\begin{eqnarray}\label{decomposition.eq}
&&F_{\tau^n_{i+1}}(_nX_{\tau^n_{i+1}-},_nA_{\tau
^n_{i+1}-})-F_{\tau
^n_i}(_nX_{\tau^n_i-},_nA_{\tau^n_i-})\nonumber\\[-2pt]
&&\qquad=F_{\tau^n_{i+1}}(_nX_{\tau
^n_{i+1}-},_nA_{\tau^n_i,h^n_i})-F_{\tau^n_i}(_nX_{\tau
^n_i},_nA_{\tau
^n_i}) \\[-2pt]
&&\qquad\quad{}+F_{\tau^n_i}(_nX_{\tau^n_{i}},_nA_{\tau^n_i-})-F_{\tau
^n_{i}}(_nX_{\tau^n_i-},_nA_{\tau^n_i-}),\nonumber
\end{eqnarray}
where we have used the fact that $F$ has predictable dependence in
the second variable to have
$F_{\tau^n_i}(_nX_{\tau^n_i},_nA_{\tau^n_i})=F_{\tau^n_i}(_nX_{\tau
^n_{i}},_nA_{\tau^n_i-})$.
The first term in \eqref{decomposition.eq} can be written $
\psi(h^n_i)-\psi(0) $ where
%
\begin{equation}
\psi(u)=F_{\tau^n_i+u}(_nX_{\tau^n_i,u},_nA_{\tau^n_i,u}).
\end{equation}
Since
$F \in\mathbb{C}^{1,2}([0,T])$, $\psi$ is right-differentiable and
left-continuous by Lem\-ma~\ref{cadlag.prop}, so:
%
\begin{eqnarray}
&&F_{\tau^n_{i+1}}(_nX_{\tau^n_i,h^n_i},_nA_{\tau^n_i,h^n_i})-F_{\tau
^n_i}(_nX_{\tau^n_i},_nA_{\tau^n_i})
\nonumber
\\[-10pt]
\\[-10pt]
\nonumber
&&\qquad=
\int_{0}^{\tau^n_{i+1}-\tau^n_i}
\mathcal{D}_{\tau^n_i+u}F(_nX_{\tau^n_i,u},_nA_{\tau^n_i,u}) \,du.
\end{eqnarray}
The second term in \eqref{decomposition.eq} can be written
$\phi(X(\tau^n_{i+1})-X(\tau^n_{i}))-\phi(0)$ where
$\phi(u)=F_{\tau^n_{i}}(_nX_{\tau^n_{i}-}^{u},_nA_{\tau^n_{i}})$.
Since $F\in\mathbb{C}^{1,2}_b$, $\phi$ is a $C^2$ function
and
$\phi'(u)=\nabla_xF_{\tau^n_{i}}(_nX_{\tau^n_{i}-}^u,_nA_{\tau
^n_{i},h_i})$, $\phi''(u)=\nabla^2_xF_{\tau^n_{i}}(_nX^u_{\tau
^n_{i}-},_nA_{\tau^n_{i},h_i})$.
Applying the It{\^o} formula to $\phi$ between 0 and
$\tau^n_{i+1}-\tau^n_i$ and the $(\mathcal{F}_{\tau_i+s})_{s \geq
0}$ continuous semimartingale $(X(\tau^n_i+s))_{s \geq0}$, yields:
%
%
\begin{eqnarray}
&&\phi\bigl(X({\tau^n_{i+1}})-X({\tau^n_{i}}) \bigr)-\phi(0)\nonumber\\[-2pt]
&&\qquad=\int_{\tau^n_{i}}^{\tau^n_{i+1}} \nabla_xF_{\tau^n_{i}}\bigl(_nX_{\tau
^n_{i}-}^{X(s)-X(\tau^n_{i})},_nA_{\tau^n_{i}} \bigr)\,dX(s) \\[-2pt]
&&\qquad\quad{}+\frac{1}{2}\int_{\tau^n_{i}}^{\tau^n_{i+1}} \operatorname{tr}\bigl[{}^t\nabla^2_xF_{\tau^n_{i}}\bigl(_nX_{\tau^n_{i}-}^{X(s)-X(\tau
^n_{i})},_nA_{\tau^n_{i}}\bigr)
\,d[X](s)\bigr].\nonumber
\end{eqnarray}
Summing over $i\geq0$ and denoting $i(s)$ the index such that $s
\in[\tau^n_{i(s)},\tau^n_{i(s)+1})$, we have shown:
%
%
\begin{eqnarray}
&&F_{t}(_nX_{t},_nA_{t})-F_{0}(X_{0},A_{0})\nonumber\\[-2pt]
&&\qquad=\int_{0}^{t}\mathcal
{D}_sF\bigl(_nX_{\tau^n_{i(s)},s-\tau^n_{i(s)}},_nA_{\tau^n_{i(s)},s-\tau
^n_{i(s)}}\bigr)\,ds
\nonumber
\\[-10pt]
\\[-10pt]
\nonumber
&&\qquad\quad{}+\int_{0}^{t} \nabla_xF_{\tau^n_{i(s)+1}}\bigl(_nX_{\tau
^n_{i(s)}-}^{X(s)-X(\tau^n_{i(s)})},_nA_{\tau^n_{i(s)},h_{i(s)}}\bigr)\,dX(s)
\\[-2pt]
&&\qquad\quad{}+ \frac{1}{2}\int_{0}^{t}\operatorname{tr}\bigl[
\nabla^2_xF_{\tau^n_{i(s)}}\bigl(_nX_{\tau^n_{i(s)}-}^{
X(s)-X(\tau^n_{i(s)})},_nA_{\tau^n_{i(s)}}\bigr)\cdot
d[X](s)\bigr]\nonumber
\end{eqnarray}
$F_{t}(_nX_{t},_nA_{t})$ converges to $F_t(X_t,A_t)$ almost surely.
Since all approximations of $(X,A)$ appearing in the various\vadjust{\goodbreak}
integrals have a $d_\infty$-distance from $(X_{s},A_{s})$ less than
$\eta_{n}\to0$, the continuity at fixed times of $\mathcal{D}F$ and
left-continuity $\nabla_xF$, $\nabla^2_{x}F$ imply that
the integrands appearing in the above integrals converge
respectively to
$\mathcal{D}_sF(X_{s},A_{s}),\nabla_x F_s(X_s,A_s),\nabla^2_{x}F_s(X_s,A_s)$
as $n\rightarrow\infty$. Since the
derivatives are in $\mathbb{B}$ the integrands in the various above
integrals are bounded by a constant dependent only on $F$, $K$ and
$R$ and $t$ does not depend on $s$ nor on $\omega$. The
dominated convergence and the dominated convergence theorem for the
stochastic integrals \cite[Chapter~IV, Theorem
32]{protter} then ensure that the Lebesgue--Stieltjes integrals converge
almost surely, and the stochastic integral in probability, to the terms
appearing in~\eqref{functional.ito.eq} as $n\rightarrow\infty$.

Consider now the general case where $X$ and $A$ may be unbounded.
Let $K_n$ be an increasing sequence of compact sets with $\bigcup_{n
\geq0} K_n = \mathbb{R}^d$ and define the optional stopping times
\[
\tau_n=\inf\{s < t| X_s \notin K^n \mbox{ or } |A_s| >
n\}\wedge t.
\]
Applying the previous result
to the stopped process $(X_{t \wedge\tau_n},A_{t\wedge\tau_n})$
and noting that, by \eqref{predictable.eq},
$F_t(X_t,A_t)=F_t(X_t,A_{t-})$ leads to
%
\begin{eqnarray*}
&&F_t(X_{t \wedge\tau
_n},A_{t\wedge\tau_n})-F_0(Z_0,A_0)\\
&&\qquad =
\int_{0}^{t\wedge\tau_n} \mathcal{D}_uF_u(X_u,A_u) \,du
+ \frac{1}{2}
\int_{0}^{t\wedge\tau_n}\operatorname{tr}({}^t\nabla^2_xF_u(X_u,A_u) \,d[X](u)) \\
&&\qquad\quad{}+ \int_{0}^{t\wedge\tau_n} {\nabla}_xF_u(X_u,A_u)\cdot dX+ \int
_{t\wedge\tau^n}^t
D_uF(X_{u \wedge\tau_n},A_{u\wedge\tau_n})\,du.
\end{eqnarray*}
The terms in the
first line converges almost surely to the integral up to time $t$
since $t\wedge\tau_n=t$ almost surely for $n$ sufficiently large.
For the same reason the last term converges almost surely to 0.
\end{pf}

\begin{remark} The above proof is probabilistic and makes use of the
(classical) It{\^o} formula~\cite{ito44}.
In the companion paper~\cite{ContFournie09c} we give a
nonprobabilistic proof of Theorem~\ref{ito.theorem}, using the
analytical approach of F\"ollmer~\cite{follmer79}, which allows $X$ to
have discontinuous (cadlag) trajectories.
\end{remark}

\begin{example} If $F_t(x_t,v_t)=f(t,x(t))$
where $f \in C^{1,2}([0,T]\times\mathbb{R}^d)$, \eqref{functional.ito.eq}~reduces to the standard It\^o formula.
\end{example}
\begin{example} For the functional in Example~\ref{doleans.ex}
$F_t(x_t,v_t)=e^{x(t)-{1}/{2}\int_0^tv(u)\,du}$,
the formula \eqref{functional.ito.eq} yields the well-known integral
representation
%
\begin{equation}
\exp\biggl(X(t)-\frac{1}{2}[X](t) \biggr)=\int_0^t
e^{X(u)-{1}/{2}[X](u)}
\,dX(u).
\end{equation}
\end{example}

An immediate corollary of Theorem~\ref{ito.theorem} is that if $X$ is a
local martingale, any $\mathbb{C}^{1,2}_b$ functional of $X$ which has
finite variation
is equal to the integral of its horizontal derivative:\vadjust{\goodbreak}
\begin{corollary} If $X$ is a local martingale and $F \in\mathbb
{C}^{1,2}_b$, the process $Y(t)=F_t(X_t,A_t)$ has finite variation if
only if $ \nabla_xF_t(X_t,A_t)=0$ $d[X]\times d\mathbb{P}$-almost everywhere.
\end{corollary}
\begin{pf} $Y(t)$ is a continuous
semimartingale by Theorem~\ref{ito.theorem}, with semimartingale decomposition
given by \eqref{functional.ito.eq}.
If $Y$ has finite variation, then by formula \eqref{functional.ito.eq},
its continuous martingale component should be zero, that is, $\int
_0^t\nabla_xF_t(X_t,A_t)\cdot dX(t)=0$ a.s.
Computing its quadratic variation, we obtain
\[
\int_0^T \operatorname{tr}({}^t\nabla_xF_t(X_t,A_t)\cdot\nabla
_xF_t(X_t,A_t)\cdot d[X])=0
\]
which implies in particular that $ \|\partial_iF_t(X_t,A_t)\|^2=0$
$d[X^i]\times d\mathbb{P}$-almost everywhere for $i=1,\ldots,d$. Thus,
$\nabla_xF_t(X_t,A_t)=0$ for $(t,\omega)\notin A\subset[0,T]\times
\Omega$ where $\int_A d[X^i]\times d\mathbb{P}=0$ for $i=1,\ldots,d$.
\end{pf}

\subsection{Vertical derivative of an adapted process}
For a ($\mathcal{F}_t$-adapted) process $Y$, the
the functional representation \eqref{Yrepresentation.eq} is not
unique, and the vertical $\nabla_xF$ depends on the choice of
representation $F$. However, Theorem~\ref{ito.theorem} implies that the
\textit{process} $\nabla_xF_t(X_t,A_t)$ has an intrinsic character, that
is, independent of the chosen representation:
\begin{corollary} \label{pathwisenabla.cor}
Let $F^1,F^2\in\mathbb{C}^{1,2}_b([0,T) )$, such that
%
\begin{equation}\forall t \in[0, T) ,\qquad F^1_t(X_t,A_t)=F^2_t(X_t,A_t)\qquad
\mathbb{P}\mbox{-a.s.}
\end{equation}
Then, outside an evanescent set,
%
\begin{eqnarray}
\qquad&&{}^{t}[\nabla_xF^1_t(X_t,A_t)-\nabla
_xF^2_t(X_t,A_t)]A(t-)[\nabla
_xF^1_t(X_t,A_t)-\nabla_xF^2_t(X_t,A_t)]
\nonumber
\\[-8pt]
\\[-8pt]
\nonumber
&&\qquad=0.
\end{eqnarray}
\end{corollary}
\begin{pf}
Let $X(t)=B(t)+M(t)$, where $B$ is a continuous process with finite
variation and $M$ is a continuous local martingale. There exists
$\Omega
_1\subset\Omega$ such that $\mathbb{P}(\Omega_1)=1$, and for
$\omega\in
\Omega$ the path of $t\mapsto X(t,\omega)$ is continuous and
$t\mapsto
A(t,\omega)$ is cadlag. Theorem~\ref{ito.theorem} implies that the
local martingale part of $0=F^1(X_t,A_t)-F^2(X_t,A_t)$ can be written as
%
\begin{eqnarray}0=\int_0^t [\nabla_xF^1_u(X_u,A_u)-\nabla_xF^2_u(X_u,A_u)
]\,dM(u).
\end{eqnarray}
Considering its quadratic variation, we have, on $\Omega_1$,
%
\begin{eqnarray}\label{pathwisenabla.eq}
0&=&\int_0^t \frac{1}{2}
{}^{t}[\nabla_xF^1_u(X_u,A_u)-\nabla_xF^2_u(X_u,A_u)]
\nonumber
\\[-8pt]
\\[-8pt]
\nonumber
&&\quad{}\times A(u-)[\nabla
_xF^1_u(X_u,A_u)-\nabla_xF^2_u(X_u,A_u)]\,du.
\end{eqnarray}
By Lemma~\ref{cadlag.prop} ($\nabla_xF^1(X_t,A_t)=\nabla
_xF^1(X_{t-},A_{t-})$ since $X$ is continuous and $F$ verifies \eqref
{predictable.eq}. So on $\Omega_1$ the integrand in \eqref
{pathwisenabla.eq} is left-continuous; therefore \eqref
{pathwisenabla.eq} implies that for $t < T$ and $\omega\in\Omega_1$,
\begin{eqnarray*}
&&{}^{t}[\nabla_xF^1_u(X_u,A_u)-\nabla_xF^2_u(X_u,A_u)]A(u-)[\nabla
_xF^1_u(X_u,A_u)-\nabla_xF^2_u(X_u,A_u)\\
&&\qquad=0.
\end{eqnarray*}
\upqed\end{pf}

In the case where for all $t<T$, $A(t-)$ is almost surely positive
definite, Corollary~\ref{pathwisenabla.cor} allows us to define
intrinsically the pathwise derivative of a process $Y$ which admits a
functional representation $Y(t)=F_t(X_t,A_t)$:

\begin{definition}[(Vertical derivative of a process)]
Define $\mathcal{C}^{1,2}_b(X)$ the set of $\mathcal{F}_t$-adapted
processes $Y$ which admit a functional representation in~$\mathbb{C}^{1,2}_b$,
%
\begin{eqnarray}\mathcal{C}^{1,2}_b(X)=\{ Y, \exists
F \in\mathbb{C}^{1,2}_b\ Y(t)=F_t(X_t,A_t)
\ \mathbb{P}\mbox{-a.s.}\} \label{Yrepresentation.eq}.
\end{eqnarray}
If $A(t)$
is nonsingular, that is, $\operatorname{det}(A(t))\neq0\ dt\times d\mathbb{P}$
almost-everywhere, then for any $Y\in\mathcal{C}^{1,2}_b(X)$, the
predictable process
\[
\nabla_XY(t)=\nabla_xF_t(X_t,A_t)
\]
is uniquely defined up to an evanescent set, independently of the
choice of $F \in\mathbb{C}^{1,2}_b$ in the representation \eqref
{Yrepresentation.eq}.
We will call $\nabla_XY$ the \textit{vertical derivative} of $Y$ with
respect to $X$.
\end{definition}

In particular this construction applies to the case where $X$ is a
standard Brownian motion, where $A=I_d$, so we obtain the existence of
a vertical derivative process for $\mathbb{C}^{1,2}_b$ Brownian functionals:
\begin{definition}[(Vertical derivative of nonanticipative Brownian functionals)]
Let $W$ be a standard d-dimensional Brownian motion. For any $Y\in
\mathcal{C}^{1,2}_b(W)$ with representation $Y(t)=F_t(W_t,t)$, the
predictable process
\[
\nabla_WY(t)=\nabla_xF_t(W_t,t)
\]
is uniquely defined up to an evanescent set, independently of the
choice of $F \in\mathbb{C}^{1,2}_b$.
\end{definition}
%

\section{Martingale representation formulas}
\label{martingaleformula.sec} Consider now the case where $X$ is a
Brownian martingale:
\begin{assumption} $ X(t)= X(0)+\int_0^t \sigma(u)\cdot dW(u)$ where
$\sigma
$ is a process adapted to $\mathcal{F}^W_t$ verifying
%
\begin{equation}
\operatorname{det}(\sigma(t))\neq0,\qquad dt\times
d\mathbb{P}\mbox{-a.e.}
\end{equation}
\end{assumption}

The functional It{\^o} formula (Theorem~\ref{ito.theorem}) then leads
to an explicit martingale representation formula for $\mathcal{F}_t$-martingales in
$\mathcal{C}^{1,2}_b(X)$. This result may be seen
as a nonanticipative counterpart of the Clark--Haussmann--Ocone\vadjust{\goodbreak}
formula~\cite{clark70,ocone84,haussmann79} and generalizes other
constructive martingale representation formulas previously obtained
using Markovian functionals
\cite{davis80,elliott88,fitz09,jacodmeleard00,pardouxpeng92},
Malliavin calculus
\cite{bismut81,karatzasocone91,haussmann79,ocone84,nualart09} or
other techniques~\cite{ahn,picard06}.

Consider an $\mathcal{F}_T$ measurable random variable $H$ with
$E|H|<\infty$,
and consider the martingale $Y(t)=E[H|\mathcal{F}_t]$.
\subsection{A martingale representation formula}
If $Y$ admits a representation $Y(t)=F_t(X_t,A_t)$ where $F \in\mathbb
{C}^{1,2}_b$, we obtain the following stochastic integral
representation for $Y$ in terms of its derivative $\nabla_XY$ with
respect to $X$:
\begin{theorem}\label{martingalerepresentation.theorem} If
$Y(t)=F_t(X_t,A_t)$ for some functional $F \in\mathbb{C}^{1,2}_b$, then
%
\begin{equation}Y(T) = Y(0) + \int_0^T \nabla_xF_t(X_t,A_t) \,dX(t)=
Y(0) + \int
_0^T \nabla_XY\cdot dX. \label{martingalerepresentation.eq}
\end{equation}
\end{theorem}

Note that regularity assumptions are not on $H=Y(T)$, but on the
martingale $Y(t)=E[H|\mathcal{F}_t],t<T$, which is typically more
regular than $H$ itself.
\begin{pf}
Theorem~\ref{ito.theorem} implies that for $t \in[0, T)$,
%
\begin{eqnarray}
Y(t)&=&\int_0^t {\mathcal D}_uF(X_u,A_u)\,du + \frac{1}{2} \int_0^t
\operatorname{tr}[{}^t\nabla^2_xF_u(X_u,A_u)\,d[X](u)]
\nonumber
\\[-8pt]
\\[-8pt]
\nonumber
&&{}+\int_0^t \nabla_xF_u(X_u,A_u)\,dX(u).
\end{eqnarray}
Given the regularity
assumptions on $F$, the first term in this sum is a continuous
process with finite variation, while the second is a continuous local
martingale. However, $Y$ is a martingale and its decomposition as
sum of a finite variation process and a local martingale is unique
\cite{revuzyor}. Hence the first term is 0, and $Y(t)=\int_0^t
F_u(X_u,A_u)\,dX_u$. Since $F\in\mathbb{C}^{0,0}_l([0,T])$ $Y(t)$ has
limit $F_T(X_T,A_T)$ as $t \rightarrow T$, so the stochastic
integral also converges.
\end{pf}

\begin{example}
If $e^{X(t)-{1}/{2}[X](t)}$ is a martingale, applying Theorem
\ref{martingalerepresentation.theorem} to the
functional $F_t(x_t,v_t)=e^{x(t)-\int_0^t v(u)\,du}$ yields the familiar formula
%
\begin{eqnarray}e^{X(t)-{1}/{2}[X](t)}=1+\int_0^t e^{X(s)-
{1}/{2}[X](s)}\,dX(s).
\end{eqnarray}
%
\end{example}
%
\subsection{Extension to square-integrable functionals}\label
{weakderivative.sec}
Let $\mathcal{L}^2(X)$ be the Hilbert space of
progressively-measurable processes $\phi$ such that
%
\begin{eqnarray}
\Vert \phi\Vert ^2_{\mathcal{L}^2(X)}=E \biggl[\int_0^t \phi_s^2 \,d[X](s)
\biggr] <\infty
\end{eqnarray}
and $\mathcal{I}^2(X)$ be the space of
square-integrable stochastic integrals with respect to $X$.
%
\begin{equation}\mathcal{I}^2(X)=\biggl\{\int_0^{.} \phi(t) \,dX(t),\phi\in
\mathcal{L}^2(X) \biggr\}
\end{equation}
endowed with the norm $\Vert Y\Vert ^2_2=E[Y(T)^2]. $
The It{\^o} integral $I_X\dvtx  \phi\mapsto\int_0^{.} \phi_s \,dX(s)$ is then
a bijective isometry from $\mathcal{L}^2(X)$ to $\mathcal{I}^2(X)$.

We will now show that the operator $\nabla_X\dvtx  \mapsto
\mathcal{L}^2(X)$ admits a suitable extension to $\mathcal{I}^2(X)$
which verifies
%
\begin{equation}\forall\phi\in\mathcal{L}^2(X),\qquad
\nabla_X\biggl( \int\phi\cdot dX \biggr) = \phi, \qquad dt \times
d\mathbb{P}\mbox{-a.s.};
\end{equation}
that is, $\nabla_X$ is the inverse of the It{\^o}
stochastic integral with respect to $X$.

%
\begin{definition}[(Space of test processes)]
The space of \textit{test processes} $D(X)$ is defined as
%
\begin{equation}
D(X)=\mathcal{C}^{1,2}_b(X) \cap
\mathcal{I}^2(X).
\end{equation}
\end{definition}

Theorem~\ref{martingalerepresentation.theorem} allows us to define
intrinsically the vertical derivative of a process in $D(X)$ as an
element of $\mathcal{L}^2(X)$.
\begin{definition}
Let $Y \in D(X)$ define the process $\nabla_XY \in\mathcal{L}^2(X)$ as
the equivalence class of $\nabla_xF_t(X_t,A_t)$, which does not depend
on the choice of the representation functional $Y(t)=F_t(X_t,A_t).$
\end{definition}

\begin{proposition}[(Integration by parts on $D(X)$)]\label{ibp.thm} Let $Y,Z \in
D(X)$. Then
%
\begin{equation}E[Y(T)Z(T)]=E\biggl[\int_0^T \nabla_XY(t)
\nabla_XZ(t)\,d[X](t) \biggr]. \label{integrationbypart.eq}
\end{equation}
\end{proposition}

\begin{pf} Let $Y,Z \in
D(X)\subset\mathcal{C}^{1,2}_b(X)$. Then $Y,Z$ are martingales with
$Y(0)=Z(0)=0$ and $E[|Y(T)|^2]<\infty, E[ |Z(T)|^2]<\infty$.
Applying Theorem~\ref{martingalerepresentation.theorem}
to $Y$ and~$Z$, we obtain
\[
E[Y(T)Z(T)]= E\biggl[ \int_0^T \nabla_XY \,dX \int_0^T\nabla
_XZ \,dX \biggr].
\]
Applying the It{\^o} isometry formula yields the result.
\end{pf}

Using this result, we can extend the
operator $\nabla_X$ to define a weak derivative on the space of (square-integrable) stochastic integrals, where $\nabla_XY$
is characterized by \eqref{integrationbypart.eq} being satisfied
against all test processes.

The following definition introduces the Hilbert space $\mathcal{W}^{1,2}(X)$
of martingales on which $\nabla_X$
acts as a weak derivative, characterized by integration-by-part formula
\eqref{integrationbypart.eq}.
This definition may be also viewed as a
nonanticipative counterpart of Wiener--Sobolev spaces in the
Malliavin calculus~\cite{malliavin,shigekawa80}.
%
\begin{definition}[(Martingale Sobolev space)]
The martingale Sobolev space $\mathcal{W}^{1,2}(X)$ is defined as
the closure in $\mathcal{I}^2(X)$
of $D(X)$.
\end{definition}

The martingale Sobolev space $\mathcal{W}^{1,2}(X)$ is in fact none
other than $\mathcal{I}^2(X)$, the set of square-integrable stochastic
integrals:
\begin{lemma} \label{density.lemma}
$\{ \nabla_XY, Y \in D(X) \}$ is dense in $\mathcal{L}^2(X)$ and
\[
\mathcal{W}^{1,2}(X)=\mathcal{I}^2(X).
\]
\end{lemma}
\begin{pf}
We first observe that the set $U$ of ``cylindrical'' processes of the form
\[
\phi_{n,f,(t_1,\ldots,t_n)}(t) = f( X(t_1),\ldots,X(t_n)) 1_{t > t_n},
\]
where $n \geq1$, $0\leq t_1<\cdots<t_n \leq T$
and $f\in C^\infty_b(\mathbb{R}^n,\mathbb{R})$ is a total set in
$\mathcal{L}^2(X)$, that is, the linear span of $U$ is dense in
$\mathcal{L}^2(X)$. For such an integrand $\phi_{n,f,(t1,\ldots,t_n)}$,
the stochastic integral with respect to $X$ is given by the
martingale
\[
Y(t)= I_X\bigl(\phi_{n,f,(t_1,\ldots,t_n)}\bigr)(t) = F_t(X_t,A_t),
\]
where the functional $F$ is defined on $\Upsilon$ as
\[
F_t(x_t,v_t)=f(x(t_1-),\ldots,x(t_n-))\bigl(x(t)-x(t_n)\bigr)1_{t > t_n},
\]
so that
\begin{eqnarray*}
\nabla_xF_t(x_t,v_t)&=&f(x_{t_1-},\ldots,x_{t_n-})1_{t > t_n},\qquad \nabla
^2_xF_t(x_t,v_t)=0,\\
 \mathcal{D}_tF(x_t,v_t)&=&0
\end{eqnarray*}
which shows that $F \in\mathbb{C}^{1,2}_b$; see Example
\ref{example.cylindrical}. Hence, $Y \in\mathcal{C}^{1,2}_b(X)$.
Since $f$ is bounded, $Y$ is obviously square integrable, so $Y\in
D(X)$. Hence $I_X(U)\subset D(X)$.

Since $I_X$ is a bijective isometry from $\mathcal{L}^2(X)$ to
$\mathcal{I}^2(X)$, the density of $U$ in $\mathcal{L}^2(X)$
entails the density of $I_X(U)$ in $\mathcal{I}^2(X)$, so ${\mathcal
W}^{1,2}(X) = \mathcal{I}^2(X)$.
\end{pf}

\begin{theorem}[(Extension of $\nabla_X$ to $\mathcal{W}^{1,2}(X)$)]
\label{weakderivative.theorem}
The vertical derivative $\nabla_X\dvtx D(X)\mapsto\mathcal{L}^2(X) $ is
closable on $\mathcal{W}^{1,2}(X)$. Its closure defines a bijective isometry
%
\begin{eqnarray}\label{closure.def}
\nabla_X\dvtx  \mathcal{W}^{1,2}(X) & \mapsto& {\mathcal
{L}}^2(X),
\nonumber
\\[-8pt]
\\[-8pt]
\nonumber
\int_0^{\cdot} \phi\cdot dX & \mapsto&
\phi
\end{eqnarray}
characterized by the following \mbox{integration by parts formula}:
for $Y \in\mathcal{W}^{1,2}(X)$, $\nabla_XY$ is the unique element of
${\mathcal{L}}^2(X)$ such that
%
\begin{equation}\label{weakderivative.eq}\quad
\forall Z \in
D(X),\qquad E[Y(T)Z(T)]=E\biggl[\int_0^T \nabla_XY(t) \nabla_XZ(t)
\,d[X](t) \biggr].
\end{equation}
In particular, $\nabla_X$ is the adjoint of the It{\^o} stochastic integral
%
\begin{eqnarray}
I_X\dvtx  \mathcal{L}^2(X) & \mapsto& \mathcal{W}^{1,2}(X),
\nonumber
\\[-8pt]
\\[-8pt]
\nonumber
\phi& \mapsto&
\int_0^{\cdot}\phi\cdot dX
\end{eqnarray}
in the following sense:
%
\begin{eqnarray}&&\forall\phi\in
\mathcal{L}^2(X), \forall Y\in\mathcal{W}^{1,2}(X),
\nonumber
\\[-8pt]
\\[-8pt]
\nonumber
&&\qquad E\biggl[ Y(T)\int_0^T\phi\cdot dX \biggr]= E\biggl[\int_0^T \nabla_XY \phi
\,d[X] \biggr].
\end{eqnarray}
\end{theorem}
\begin{pf} Any $Y \in\mathcal{W}^{1,2}(X)$
may be written as $ Y(t)=\int_0^t \phi(s) \,dX(s) $ with $\phi\in
\mathcal{L}^2(X)$, which is uniquely defined $d[X]\times
d\mathbb{P}$ a.e. The It{\^o} isometry formula then guarantees that
\eqref{weakderivative.eq} holds for $\phi$. To show
that \eqref{weakderivative.eq} uniquely characterizes $\phi$, consider
$\psi\in\mathcal{L}^2(X)$ which also satisfies \eqref
{weakderivative.eq}, then, denoting
$I_X(\psi)=\int_0^{\cdot}\psi \,dX$ its stochastic integral with respect to
$X$, \eqref{weakderivative.eq} then implies that
\[
\forall Z\in D(X), \qquad \langle I_X(\psi)-Y,Z\rangle _{\mathcal{W}^{1,2}(X)}=E\biggl[
\biggl(Y(T)-\int_0^T\psi \,dX\biggr)Z(T) \biggr]=0
\]
which implies $I_X(\psi)= Y$ $d[X]\times d\mathbb{P}$ a.e., since by
construction $D(X)$ is dense in $\mathcal{W}^{1,2}(X)$.
Hence, $\nabla_X\dvtx  D(X) \mapsto\mathcal{L}^2(X)$ is closable on $\mathcal
{W}^{1,2}(X)$.

This
construction shows that $\nabla_X\dvtx  \mathcal{W}^{1,2}(X)\mapsto
{\mathcal
{L}}^2(X)$ is a bijective
isometry which coincides with the adjoint of the
It{\^o} integral on $\mathcal{W}^{1,2}(X)$.
\end{pf}

Thus, the It{\^o} integral $I_X$ with respect to
$X$,
\[
I_X\dvtx  \mathcal{L}^2(X)  \mapsto \mathcal{W}^{1,2}(X),
\]
admits an inverse on $\mathcal{W}^{1,2}(X)$ which is
an extension of the (pathwise) vertical derivative $\nabla_X$
operator introduced in Definition~\ref{verticalderivative.def}, and
%
\begin{equation}
\forall\phi\in\mathcal{L}^2(X),\qquad \nabla_X \biggl(
\int_0^{\cdot}\phi \,dX \biggr)=\phi
\end{equation}
holds in the sense of equality in
$\mathcal{L}^2(X)$.

The above results now allow us to state a general version of the
martingale representation formula, valid for all square-integrable
martingales:
\begin{theorem}[(Martingale representation formula: general
case)]\label{L2martingalerepresentation.theorem} For any
square-integrable $(\mathcal{F}_t^X)_{t\in[0,T]}$-martingale $Y$,
\[
Y(T) = Y(0) + \int_0^T \nabla_XY \,dX,\qquad \mathbb{P}\mbox{-}a.s.
\]
\end{theorem}
%
\section{Relation with the Malliavin
derivative}\label{sec.malliavin} The above results hold in
particular in the case where $X=W$ is a Brownian motion. In this
case, the vertical derivative $\nabla_W$ may be related to the
\textit{Malliavin derivative}~\cite{malliavin,bismut81,bismut83,stroock81}
as follows.

Consider the canonical Wiener space
$(\Omega_0=C_0([0,T],\mathbb{R}^d),\|\cdot\|_\infty,\mathbb{P})$ endowed
with its Borelian $\sigma$-algebra, the filtration of the canonical
process. Consider an $\mathcal{F}_T$-measurable functional
$H=H(X(t),t\in[0,T])=H(X_T)$ with $E[|H|^2]<\infty$. If $H$ is
differentiable in the Malliavin sense
\cite{bismut81,malliavin,nualart09,stroock81}, for example, $H\in{\bf
D}^{1,2}$ with Malliavin derivative $\mathbb{D}_tH$, then the
Clark--Haussmann--Ocone formula~\cite{ocone84,nualart09} gives a
stochastic integral representation of $H$ in terms of the Malliavin
derivative of $H$.
%
\begin{equation}H= E[H]+\int_0^T {}^p\hspace*{-1.5pt}E[ \mathbb{D}_tH
|\mathcal{ F}_t] \,dW_t,\label{clarkocone.eq}
\end{equation}
where ${}^p\hspace*{-1.5pt}E[
\mathbb{D}_tH |\mathcal{ F}_t]$ denotes the predictable projection
of the Malliavin derivative. This yields a stochastic integral
representation of the martingale $Y(t)=E[H|\mathcal{F}_t]$,
\[
Y(t)=E[H|\mathcal{F}_t]= E[H]+\int_0^t {}^p\hspace*{-1.5pt}E[ \mathbb{D}_tH
|\mathcal{ F}_u] \,dW_u.
\]
Related martingale representations have
been obtained under a variety of conditions
\cite{bismut81,davis80,fitz09,karatzasocone91,pardouxpeng92,nualart09}.

Denote by:
\begin{itemize}
\item$L^2([0,T]\times\Omega)$ the set of (anticipative)
processes $\phi$ on $[0,T]$ with\break $E\int_0^T\|\phi(t)\|^2 \,dt<\infty$;
\item
$\mathbb{D}$ the
Malliavin derivative operator, which associates to a random
variable $H \in\mathbf{D}^{1,2}(0,T)$ the (anticipative) process
$(\mathbb{D}_tH)_{t\in[0,T]}\in L^2([0,T]\times\Omega)$.
\end{itemize}

\begin{theorem}[(Lifting theorem)]\label{relevement.theorem}
The following diagram is commutative is the sense of $dt\times
d\mathbb{P}$ equality:
\begin{center}$
\begin{array}[c]{ccc}
\mathcal{I}^2(W)&\stackrel{\nabla_W}{\rightarrow}&\mathcal
{L}^2(W)\\
\uparrow\scriptstyle{(E[\cdot|\mathcal{F}_t])_{t \in[0, T]}}&&\uparrow
\scriptstyle{(E[\cdot|\mathcal{F}_t])_{t\in[0, T]}}\\
\mathbf{D}^{1,2}&\stackrel{\mathbb{D}}{\rightarrow}&
L^2([0,T]\times\Omega
).\\
\end{array}
$
\end{center}
In other words, the conditional expectation operator intertwines
$\nabla_W$ with the Malliavin derivative,
%
\begin{eqnarray}\forall H \in
L^2(\Omega_0,\mathcal{F}_T,\mathbb{P}),\qquad
\nabla_W(E[H|\mathcal{F}_t]) =
E[\mathbb{D}_tH|\mathcal{F}_t].
\end{eqnarray}
\end{theorem}
\begin{pf}
The Clark--Haussmann--Ocone formula~\cite{ocone84} gives
%
\begin{eqnarray}\forall H
\in\mathbf{D}^{1,2}, \qquad H= E[H]+\int_0^T {}^p\hspace*{-1.5pt}E[
\mathbb{D}_tH |\mathcal{ F}_t] \,dW_t,
\end{eqnarray}
where ${}^p\hspace*{-1.5pt}E[ \mathbb{D}_tH
|\mathcal{ F}_t]$ denotes the predictable projection of the
Malliavin derivative. On other hand, Theorem
\ref{martingalerepresentation.theorem} gives
%
\begin{eqnarray}\forall H \in
{L}^{2}(\Omega_0,\mathcal{F}_T,\mathbb{P}),\qquad  H= E[H] + \int_0^T
\nabla_W Y(t) \,dW(t),
\end{eqnarray}
where $Y(t)=E[H|\mathcal{F}_t]$. Hence $
{}^p\hspace*{-1.5pt}E[ \mathbb{D}_tH
|\mathcal{ F}_t]=\nabla_WE[H|\mathcal{F}_t] $, $dt \times
d\mathbb{P}$ almost everywhere.
\end{pf}

Thus, the conditional expectation operator (more precisely: the
\textit{predictable} projection on $\mathcal{F}_t$ \cite[Volume I]{dm})
can be viewed as a morphism which ``lifts''
relations obtained in the
framework of Malliavin calculus into relations between
nonanticipative quantities, where the Malliavin derivative and the
Skorokhod integral are replaced, respectively, by the vertical derivative
$\nabla_W$ and the It{\^o} stochastic integral.

From a computational viewpoint, unlike the Clark--Haussmann--Ocone
representation which requires to simulate the \textit{anticipative}
process $\mathbb{D}_tH$ and compute conditional expectations,
$\nabla_X Y$ only involves nonanticipative quantities which can be
computed path by path. It is thus more amenable to numerical
computations. This topic is further explored in a forthcoming work.

\begin{appendix}
\section*{Appendix: Proof of Theorem 2.7}
\label{measurabilityproof.sec}
In order to prove Theorem~\ref{Measurability} in the general case where
$A$ is only required to be cadlag, we need the following three
lemmas. The first lemma states a property analogous to ``uniform
continuity'' for cadlag functions: 
%
\renewcommand{\thelemmaa}{A.\arabic{lemmaa}}
\setcounter{lemmaa}{0}
\begin{lemmaa}\label{uniform.lemma} 
Let $f$ be a cadlag function on $[0,T]$ and define $\Delta
f(t)=f(t)-f(t-)$. Then
%
%
\begin{eqnarray}\label{uniformcadlag.eq}
&&\forall\varepsilon> 0, \exists\eta(\varepsilon) >0,
\nonumber
\\[-8pt]
\\[-8pt]
\nonumber
&&\qquad |x-y|\leq
\eta\quad\Rightarrow\quad|f(x)-f(y)|\leq\varepsilon+\sup_{t\in(x,y]}\{|\Delta
f(t)|\}.
\end{eqnarray}
\end{lemmaa}
\begin{pf} If \eqref{uniformcadlag.eq} does not hold, then there
exists a sequence $(x_{n},y_{n})_{n\geq1}$ such that $x_{n}\leq
y_{n}$, $y_{n}-x_{n}\rightarrow0$, but $|f(x_{n})-f(y_{n})|>\varepsilon
+\sup_{t\in[x_{n},y_{n}]}\{|\Delta f(t)|\}$. We can extract a
convergent subsequence $(x_{\psi(n)})$ such that $x_{\psi
(n)}\rightarrow x$. Noting that either an infinity of terms of the
sequence are less than $x$ or an infinity are more than $x$, we can
extract \textit{monotone} subsequences $(u_{n},v_{n})_{n\geq1}$ which
converge to $x$.
If $(u_n),(v_n)$ both converge to $x$ from above or from below,\vadjust{\goodbreak}
$|f(u_{n})-f(v_{n})|\rightarrow0$ which yields a contradiction. If one
converges from above and the other from below,
$\sup_{t\in[u_{n},v_{n}]}\{|\Delta f(t)|\}\geq|\Delta f(x)|$, but
$|f(u_{n})-f(v_{n})|\rightarrow|\Delta f(x)|$, which results in a
contradiction as well. Therefore \eqref{uniformcadlag.eq} must hold.
\end{pf}
\begin{lemma}\label{stoppingtime.lemma}
If $\alpha\in\mathbb{R}$ and $V$ is an adapted cadlag process defined
on a filtered probability space $(\Omega,\mathcal{F},(\mathcal
{F}_{t})_{t\geq0},\mathbb{P})$ and $\sigma$ is a optional time, then
%
%
\begin{equation}
\tau=\inf\{t > \sigma, |V(t)-V(t-)|> \alpha\}
\end{equation}
is a stopping time.
\end{lemma}

\begin{pf}
We can write that
%
%
\begin{equation}
\{\tau\leq t\}=\bigcup_{q\in\mathbb{Q}\cap[0,t)}(\{\sigma\leq
t-q\}\cap\Bigl\{\sup_{t\in(t-q,t]}|V(u)-V(u-)|> \alpha\Bigr\}
\end{equation}
and, using Lemma~\ref{uniform.lemma},
%
%
\begin{eqnarray}
&&\Bigl\{\sup_{u\in(t-q,t]}|V(u)-V(u-)|> \alpha\Bigr\}
\nonumber
\\[-4pt]
\\[-12pt]
\nonumber
&&\qquad=\bigcup_{n_0>1}
\bigcap_{n>n_0} \bigcup_{m\geq1}\biggl\{\sup_{1\leq i\leq
2^n}\biggl|V\biggl(t-q\frac{i-1}{2^n}\biggr)-V\biggl(t-q\frac{i}{2^n}\biggr)\biggr|> \alpha+
\frac{1}{m}\biggr\}.
\end{eqnarray}
\upqed\end{pf}
\begin{lemma}[(Uniform approximation of cadlag functions by step
functions)]\label{Approximation}
Let $f\in D([0,T],\mathbb{R}^d)$ and $\pi^n=(t^n_i)_{n\geq1,
i=0,\ldots,k_n}$ a sequence of partitions
$(0=t^{n}_{0}<t_{1}<\cdots<t^{n}_{k_{n}}=T)$ of $[0,T]$ such that
\[
\mathop{\sup}_{0 \leq i \leq k_n-1} |t^{n}_{i+1}-t^{n}_{i}| \mathop
{\to
}^{n\rightarrow\infty} 0,\qquad \mathop{\sup}_{u \in[0,T]\setminus\pi
^n} |\Delta f(u) | \mathop{\to}^{n \rightarrow\infty} 0
\]
then
\begin{equation}
\sup_{u \in[0,T]} \Biggl|f(u)-\sum_{i=0}^{k_{n}-1}
f(t^n_{i})1_{[t^{n}_{i},t^{n}_{i+1})}(u)+f(t^{n}_{k_{n}})1_{\{
t^{n}_{k_{n}}\}}(u)\Biggr|
\mathop{\to}^{n\rightarrow\infty} 0.
\end{equation}
\end{lemma}
\begin{pf} Denote $h^n=f-\sum_{i=0}^{k_{n}-1}
f(t^n_{i})1_{[t^{n}_{i},t^{n}_{i+1})}+f(t^{n}_{k_{n}})1_{\{
t^{n}_{k_{n}}\}}$.
Since $f-h^n$ is piecewise constant on $\pi^n$ and $h^n(t_i^n)=0$ by
definition,
\[
\sup_{t\in[0,T]}|h^n(t)|=\sup_{i=0,\ldots,k_n-1}\sup
_{[t_i^n,t_i^{n+1})}|h^n(t)|=\sup_{t_i^n<t<t_i^{n+1}}|f(t)-f(t_i^n)|.
\]
Let $\varepsilon>0$. For $n\geq N$ sufficiently large, $\sup_{u \in
[0,T]\setminus\pi^n } |\Delta f(u) |\leq\varepsilon/2$ and\break $\sup
_i|t^{n}_{i+1}-t^{n}_{i}|\leq\eta(\varepsilon/2)$ using the notation of
Lemma~\ref{uniform.lemma}.
Then, applying Lemma~\ref{uniform.lemma} to $f$ we obtain, for $n\geq N$,
\[
\sup_{t\in[t_i^n,t_i^{n+1})}|f(t)-f(t_i^n)|\leq\frac{\varepsilon}{2} +
\sup_{t_i^n<t<t_i^{n+1}} |\Delta f(u) |\leq\varepsilon.
\]
\upqed\end{pf}

We can now prove Theorem~\ref{Measurability} in the case where $A$
is a cadlag adapted process.\vadjust{\goodbreak}

\begin{pf*}{Proof of Theorem \protect\ref{Measurability}} Let us first show that
$F_t(X_t,A_t)$ is adapted. Define
%
%
\begin{equation}
\qquad\tau^N_{0}=0,\qquad
\tau^N_{k}=\inf\biggl\{t > \tau^N_{k-1} | 2^{N}t \in\mathbb{N} \mbox{ or } |A(t)-A(t-)| > \frac{1}{N} \biggr\}\wedge t.
\end{equation}
From lemma~\ref{stoppingtime.lemma}, $\tau^N_{k}$ are stopping times.
Define the following piecewise constant approximations of $X_{t}$ and
$A_{t}$ along the partition $(\tau^N_{k},k\geq0)$:
%
%
\begin{eqnarray}
X^{N}(s)&=&\sum_{k\geq0} X_{\tau^N_{k}}1_{[\tau^N_{k},\tau
^N_{k+1}[}(s)+X(t)1_{\{t\}}(s),
\nonumber
\\[-8pt]
\\[-8pt]
\nonumber
A^{N}(s)&=&\sum_{k=0} A_{\tau^N_{k}}1_{[\tau^N_{k},\tau
^N_{k+1})}(t)+A(t)1_{\{t\}}(s),
\end{eqnarray}
as well as their truncations of rank $K$
%
%
\begin{equation}\quad
_{K}X^{N}(s)=\sum_{k=0}^{K} X_{\tau^N_{k}}1_{[\tau^N_{k},\tau
^N_{k+1})}(s),\qquad
{}_{K}A^{N}(t)=\sum_{k=0}^{K} A_{\tau^N_{k}}1_{[\tau^N_{k},\tau
^N_{k+1})}(t).
\end{equation}
Since $(_{K}X^{N}_{t},_{K}A^{N}_{t})$ coincides with
$(X^{N}_{t},A^{N}_{t})$ for $K$ sufficiently large,
%
%
\begin{equation}
F_t(X^{N}_{t},A^{N}_{t})=\lim_{K\rightarrow\infty}
F_{t}(_KX^{N}_{t},_KA^{N}_{t}).
\end{equation}
The approximations
$F^n_{t}(_{K}X^{N}_{t},_{K}A^{N}_{t}) $
are $\mathcal{F}_t$-measurable as they are
continuous functions of the random variables
\[
\{(X(\tau^N_{k})1_{\tau^N_{k}\leq
t},A(\tau^N_{k})1_{\tau^N_{k}\leq t}),k \leq K\},
\]
so their limit $F_{t}(X^{N}_{t},A^{N}_{t})$ is also $\mathcal
{F}_t$-measurable. Thanks to Lemma~\ref{Approximation},
$X^{N}_{t}$ and $A^{N}_{t}$ converge uniformly to $X_{t}$ and
$A_{t}$, and hence $F_{t}(X^{N}_{t},A^{N}_{t})$ converges to
$F_t(X_{t},A_{t})$ since $F_t\dvtx (D([0,t],\mathbb{R}^d)\times\mathcal
{S}_t, \|\cdot\|_\infty)\to\mathbb{R}$ is continuous.

To show the optionality of $Z$ in point (ii), we will show that $Z$ it
as limit of right-continuous adapted processes. For $t \in[0,T]$,
define $i^n(t)$ to be the integer such that $t \in[\frac{iT}{n},\frac
{(i+1)T}{n}).$ Define the process
$Z^n_t=F_{{(i^n(t))T}/{n}}(X_{{(i^n(t))T}/{n}},\break
A_{{(i^n(t))T}/{n}})$, which is piecewise-constant and has right-continuous
trajectories, and is also adapted by the first part of the theorem.
Since $F\in\mathbb{C}^{0,0}_l$, $Z^n(t) \rightarrow Z(t)$ almost
surely, which proves that $Z$ is optional.
Point (iii) follows from (i) and Lemma~\ref{cadlag.prop}, since in both
cases $F_t(X_t,A_t)=F_t(X_{t-},A_{t-})$, and hence $Z$ has
left-continuous trajectories.
\end{pf*}
\end{appendix}

\section*{Acknowledgments}
We thank Bruno Dupire for sharing his original ideas with
us, Hans-J\"urgen Engelbert, Hans F\"ollmer, Jean Jacod, Shigeo Kusuoka
and an anonymous referee for helpful comments.
R. Cont is especially grateful to the late Paul Malliavin for
encouraging this work.

%

\printaddresses

\end{document}